\documentclass[10pt]{article}

\usepackage{a4,fancyhdr}
\usepackage{amsmath,amssymb,latexsym}
\usepackage{array,graphicx}
\usepackage{hyperref,stmaryrd,xcolor} 
\usepackage{pgf,tikz}
\usetikzlibrary{arrows}

\newtheorem{theorem}{Theorem}

\newtheorem{lemma}{Lemma}

\newcommand{\va}{\varphi}
\newcommand{\ppp}{\partial}
\newcommand{\www}{\widetilde}

\def\endproof{\hfill{}$\square$}

\renewcommand{\thefootnote}{\fnsymbol{footnote}}

\author{Lucie Baudouin\footnotemark[1]\,\,, Masahiro Yamamoto\footnotemark[2]}

\begin{document}

\title{ \bf Inverse problem on a tree-shaped network}

\maketitle

\footnotetext[1]{CNRS ; LAAS ; 7 avenue du colonel Roche, F-31077 Toulouse, France ;\\
Universit\'e de Toulouse ; UPS, INSA, INP, ISAE, UT1, UTM, LAAS ; F-31077 Toulouse, France.\\
 E-mail: {\tt lucie.baudouin@laas.fr}}
\footnotetext[2]{Graduate School of Mathematical Sciences, University of Tokyo, 3-8-1 Komaba, Tokyo, 153-8914 Japan.\\
 E-mail: {\tt myama@next.odn.ne.jp}}

\abstract{
In this article, we prove a uniqueness result for a coefficient inverse problems regarding a wave, a heat 
or a Schr\"odinger equation set on a tree-shaped network, as well as the corresponding stability result of the inverse problem for the wave equation.
The objective is the determination of the potential on each edge of the network from the additional
measurement of the solution at all but one external end-points.
Our idea for proving the uniqueness is to use a traditional approach 
in coefficient inverse problem by Carleman estimate.
Afterwards, using an observability estimate on the whole 
network, we apply a compactness-uniqueness argument and prove 
the stability for the wave inverse problem. 
}
\vspace{0.2 cm}


\vspace{0.2 cm} {\bf AMS subject classifications:} 35R30, 93C20, 34B45

\renewcommand{\thefootnote}{\fnsymbol{footnote}}

\section{Introduction and main results}

Let $\Lambda$ be a tree-shaped network composed of
$ N+1$ open segments $(e_j)_{j=0, 1, ..., N}$ of length~$\ell_j$, linked 
by $N_1$ internal node points belonging to the set $\Pi_1$ and let us denote 
by $\Pi_2$ the set of  $N_2$ exterior end-points where only one segment starts.
Here we note that $N+1 = N_1+N_2$.
By ``tree-shaped network'', we mean that $\Lambda$ does not contain any closed 
loop.

For any function $f: \Lambda\rightarrow\mathbb{R}$ and any internal node 
$P\in\Pi_1$ where $n_P$ segments, say $e_1, ..., e_{n_P}$, meet, we set
$$
f_j=f\vert_{e_j}:  \mbox{the restriction of $f$ to the edge $e_j$}, \mbox{ and}
\quad
\left[ f \right]_P := \sum_{j=1}^{n_P}f_{j}(P).
$$
We consider on this plane $1$-d tree-shaped network $\Lambda$ either wave or 
heat or even Schr\"odinger equations, with a different potential term $x\mapsto p_j(x)$ on each 
segment.

Our first, and main, system of interest is the following $1$-d wave equation on the network $\Lambda$:
\begin{equation}\label{NW}
\left\{
\begin{array}{lll}
\ppp_t^2 u_j - \ppp_x^2 u_j + p_j(x)u_j = 0\quad &\quad\forall 
j\in\{0,1,..., N\}, (x,t)\in e_j \times (0,T),\\
u(Q,t)=h(t),&\quad \,\forall Q\in \Pi_2, t\in(0,T),\\
u(x,0) = u^0(x), \ppp_t u(x,0)=u^1(x),&\quad x\in\Lambda,
\end{array}
\right.
\end{equation}
assuming some compatibility condition between the boundary and initial data.
Moreover we assume the continuity and what is called the Kirchhoff law at 
any internal node $P\in\Pi_1$, which are given by
\begin{equation}\label{C}
u_j(P,t)=u_i(P,t)=:u(P,t),\quad\forall i,j\in\left\{1,...,n_P\right\},\,
0<t<T,\\
\end{equation}
\begin{equation}\label{K}
\left[ u_x (t)\right]_P :=\sum_{j=1}^{n_P} \ppp_x u_{j}(P,t)=0,\quad 0<t<T.
\end{equation}
Henceforth we choose an orientation of $\Lambda$ such that
to two endpoints of each segment $e$, correspond an initial node $I(e)$ and 
a terminal
node $T(e)$. 
We further define the outward normal derivative 
$\ppp_{n_e}u_j$ at a node $P$ of $e_j$ by 
$$
\ppp_{n_e}u_j(P,t) =
\left\{
\begin{array}{lll}
-\ppp_xu_j(P,T), \quad &\mbox{if $P\in I(e_j)$}, \\
\ppp_xu_j(P,T), \quad &\mbox{if $P\in T(e_j)$}. \\
\end{array}
\right.
$$
Henceforth we set 
$$
u = (u_0, ..., u_N), \quad u_j = u\vert_{e_j}, \quad\hbox{ and } \quad
p = (p_0, ..., p_N), \quad p_j = p\vert_{e_j} \quad
\mbox{for $j \in \{0, 1, ..., N\}$}.
$$
Let us also mention that at a node point, at least three segments $e_j$ meet.
If only two segments, say $e_1, e_2$, meet at a node point, then by \eqref{C} 
and \eqref{K}, setting
$u = u_1$ and $p=p_1$ in $e_1$ and $u = u_2$, $p=p_2$ in $e_2$, we have
$\ppp_t^2u - \ppp_x^2 u + p u$ in $e_1 \cup e_2$.  Therefore we can regard
$e_1 \cup e_2$ as one open segment.\\
Since one can prove the unique existence of solution to (1) - (3) 
in a suitable function space (e.g., Lions and Magenes \cite{LionsMagenesBook}), we denote the solution by $u[p](x,t)$, 
and we set $u[p] = (u[p]_0, ..., u[p]_N)$.
 
Moreover we consider the following heat system on the same network~$\Lambda$ 
\begin{equation}\label{NH}
\left\{
\begin{array}{lll}
\ppp_t u_j - \ppp_x^2 u_j + p_j(x)u_j = 0\quad &\quad\forall j\in\{0,1,..., 
N\}, \forall (x,t)\in e_j \times (0,T),\\
\ppp_xu(Q,t)=0,&\quad \forall Q\in \Pi_2, \forall t\in(0,T),\\
u(x,0) = u^0(x), &\quad \forall x\in\Lambda,
\end{array}
\right.
\end{equation}
and the Schr\"odinger system on the network $\Lambda$ 
\begin{equation}\label{NS}
\left\{
\begin{array}{lll}
i\ppp_t u_j - \ppp_x^2 u_j + p_j(x)u_j = 0\quad &\quad\forall 
j\in\{0,1,..., N\}, \forall (x,t)\in e_j \times (0,T),\\
u(Q,t)=h(t),&\quad \forall Q\in \Pi_2, \forall t\in(0,T),\\
u(x,0) = u^0(x), &\quad \forall x\in\Lambda,
\end{array}
\right.
\end{equation}
both under the same node conditions \eqref{C} and \eqref{K}.
Here and henceforth we set $i = \sqrt{-1}$.
If there is no possible confusion, by the same notation 
$u[p]$ we denote the solution to \eqref{NH} or \eqref{NS}, under \eqref{C} and \eqref{K}.\\

\noindent\textbf{Inverse Problem}: Is it possible to retrieve the potential 
$p$ everywhere in the whole network~$\Lambda$ from measurements at all 
external nodes except one?\\

In our article, we address the following two fundamental theoretical questions
concerning coefficient inverse problems: \\

\noindent\textbf{Uniqueness}: Do the equalities of the measurements 
$\ppp_x u[p] (Q,t) = \ppp_x u[q] (Q,t)$ for all $t\in(0,T)$ and $Q\in
\Pi_2\setminus {\{Q_{N_2}\}}$ imply $p = q$ on $\Lambda$?\\

\noindent\textbf{Stability}: Can we estimate, in appropriate norms, the 
difference of two potentials $p - q$ on $\Lambda$ by the difference of the 
corresponding measurements $\ppp_x u[p] (Q,t) - \ppp_x u[q] (Q,t)$  for all 
$t\in(0,T)$ and $Q\in \Pi_2\setminus {\{Q_{N_2}\}}$ ?\\

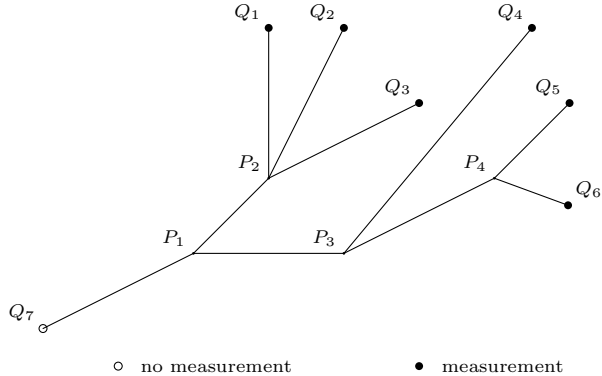
\begin{figure}[ht]
    \begin{center}
\definecolor{qqqqff}{rgb}{0,0,1}
\begin{tikzpicture}[line cap=round,line join=round,>=triangle 45,x=1.0cm,y=1.
0cm]
\clip(-1,-1) rectangle (8,5);
\draw (0,0)-- (2,1);
\draw (2,1)-- (3,2);
\draw (3,2)-- (3,4);
\draw (3,2)-- (4,4);
\draw (3,2)-- (5,3);
\draw (2,1)-- (4,1);
\draw (4,1)-- (6.5,4);
\draw (4,1)-- (6,2);
\draw (6,2)-- (7,3);
\draw (6,2)-- (6.98,1.64);
\begin{scriptsize}
\draw (0,0) circle (1.5pt);
\draw (0,0) node[above left] {$Q_7$};
\fill (2,1) circle (0.5pt);
\draw (2,1) node[above left] {$P_1$};
\fill  (3,2) circle (0.5pt);
\draw (3,2) node[above left] {$P_2$};
\fill  (3,4) circle (1.5pt);
\draw (3,4) node[above left] {$Q_1$};
\fill  (4,4) circle (1.5pt);
\draw (4,4) node[above left] {$Q_2$};
\fill  (5,3) circle (1.5pt);
\draw (5,3) node[above left] {$Q_3$};
\fill  (4,1) circle (0.5pt);
\draw (4,1) node[above left] {$P_3$};
\fill  (6.5,4) circle (1.5pt);
\draw (6.5,4) node[above left] {$Q_4$};
\fill  (6,2) circle (0.5pt);
\draw (6,2) node[above left] {$P_4$};
\fill  (7,3) circle (1.5pt);
\draw (7,3) node[above left] {$Q_5$};
\fill  (6.98,1.64) circle (1.5pt);
\draw  (6.98,1.64)  node[above right] {$Q_6$};
\draw  (1,-0.5) circle (1.5pt);
\draw (1.2,-0.5) node[right] {no measurement};
\fill  (5,-0.5) circle (1.5pt);
\draw (5.2,-0.5) node[right] {measurement};
\end{scriptsize}
\end{tikzpicture}
\caption{A star-shaped network with $10$ edges ($N=9$, $N_1 = 4$, $N_2= 7$).
        \label{fig}}
    \end{center}
\end{figure}

This inverse problem is nonlinear and we will  give here the proof of 
the uniqueness of the solution with an argument which do not use a global 
Carleman estimate.
Very recent papers on coefficient inverse problems on networks, as
Baudouin, Cr\'epeau and Valein \cite{BaudouinCrepeauValein11} for the wave 
equation, and Ignat, Pazoto and Rosier \cite{IgnatPazRosier12} for the heat and the Schr\"odinger equations, give 
indeed \textit{stability} and therefore  \textit{uniqueness} from appropriate global 
Carleman estimates. Our first goal is to prove 
the uniqueness of the potential on the tree-shaped network 
from measurements only at all the exterior end-points of the network, except one.  
The argument for the uniqueness  
will work for either the wave or the heat or the Schr\"odinger equations on 
the network.  The question of the proof of the Lipschitz stability in the case 
of the wave equation will be addressed afterwards, using a compacteness-uniqueness argument, and relies on the 
observability estimate on the whole network which was already proved in the 
literature in several situations. 

Concerning the precise topic which we are considering, the bibliography lies in
two different domains, namely coefficient inverse problems for partial differential equation on the one hand and control and stabilization in networks on the 
other hand.

Therefore one can begin by mentioning the book of Isakov \cite{Isakov} which 
adresses some techniques linked to the study of inverse problem for several 
partial differential equations.
Actually, as the first answer to the uniqueness for a coefficient inverse 
problem with a single measurement, we refer to Bukhgeim and Klibanov 
\cite{BuKli81}, and see also 
Klibanov \cite{Klibanov92} and Yamamoto \cite{Yam99} for example.  
Here we  do not intend to give an exhaustive list of references.
After the proof of uniqueness using the basic $1$-d result on the basis of
local Carleman estimates, the idea beneath this article is to take advantage 
of an observability estimate to obtain the Lipschitz stability of the inverse 
problem with a compactness-uniqueness argument.
Nowadays, many results on the stability of inverse problems are derived 
directly from global Carleman estimates, and see e.g., 
\cite{BaudouinCrepeauValein11} and \cite{IgnatPazRosier12}. 
One should also know that studies on inverse problems and controllability of 
partial differential equations share some technical materials such as 
Carleman estimates and observability inequalities.
In the particular network setting, we would like to make use of 
classical results such as well-known $1$-d local Carleman 
estimates, observability estimates on the network borrowed from control 
studies, in order to obtain uniqueness and stability results.
We can also give some more references on inverse problems for hyperbolic 
equations such as Baudouin, Mercado and Osses \cite{BMO07}, 
Imanuvilov and Yamamoto \cite{ImYamIP01}, \cite{ImYamCom01}, Puel and Yamamoto
\cite{PuelYam97}, Yamamoto and Zhang \cite{YamZhang03}, which are all based 
upon local or global Carleman estimates. 

Besides, the control, observation and stabilization problems of networks have 
been the object of recent and intensive researches such as e.g., 
D\'ager and Zuazua \cite{DagZua06}, Lagnese, Leugering and Schmidt 
\cite{LaLeSch94}, Zuazua \cite{ZuaSurveyNetworks}. 
More specifically, the control being only applied at one single end 
of the network, the articles D\'ager \cite{Dag04}, D\'ager and Zuazua 
\cite{DagZua00,DagZua06} 
prove controllability results for the wave equation on networks, using 
observability inequalities under assumptions about the irrationality properties
 of the ratios of the lengths of the strings. We can also underline that many 
results of controllability on networks concern only the wave equation without 
lower order terms (see \cite{LaLeSch94}, Schmidt \cite{Schmidt92} 
for instance). 
However it is difficult to consider such measurements at more limited 
nodes for the inverse problem and we do not consider the measurements
at less external nodes. 
\\

In the sequel, we shall use the following notations: 
\begin{eqnarray*}
L^{\gamma}(\Lambda)&=&\left\{f; \thinspace f_j\in L^{\gamma}(e_j),
\,\forall j\in\{0,1,..., N\}\right\}, \quad \gamma \ge 1,\\
H^1_0(\Lambda)&=& \Big\{f; \thinspace f_j\in 
H^1(e_j),\,\forall j\in\{0,1,..., N\}, \, f_j(P)=f_k(P) 
\thinspace \mbox{if $e_j$ and $e_k$ meet at $P$},\\
		&&\, \forall P \in\Pi_1,\, 
\textnormal{and } f(Q)=0, \, \forall Q\in \Pi_2 \Big\}.
\end{eqnarray*}
For shortness, for $f\in L^1(\Lambda)$, we often write,
\begin{equation*}
\int_{\Lambda}f dx=\sum_{j=0}^N \int_{e_j} f_j(x)dx,
\end{equation*}
where the integral on $e_j$ is oriented from $I(e_j)$ to $T(e_j)$.
Then the norms of the Hilbert spaces $L^2(\Lambda)$ and $H_0^1(\Lambda)$ are 
defined by
$$
\left\|f\right\|_{L^2(\Lambda)}^2=\int_{\Lambda}\left|f\right|^2dx
\hbox{ and }\left\|f\right\|_{H_0^1(\Lambda)}^2=\int_{\Lambda}\left|\ppp_x
f\right|^2dx.
$$
For $M\ge 0$, we introduce the set
$$
L^\infty_M(\Lambda) =  \left\{q=(q_0,..., q_N); \thinspace q_j\in 
L^\infty(e_j),\, 
\forall j\in\{0,1,..., N\}\thinspace 
\mbox{such that $\|q\|_{L^\infty(\Lambda)} 
\leq M$} \right\}.
$$

We are ready to state our first main result:

\begin{theorem}[\bf Uniqueness]\label{Thm1}
Let $r>0$ be an arbitrary constant.  
Assume that $p, q \in L^{\infty}(\Lambda)$ and the initial value $u^0$ 
satisfies 
$$ 
	|u^0(x)|\geq r> 0, \quad \mbox{a.e. in $\Lambda$}.
$$
Assume further that the solutions $u[p], u[q]$ of \eqref{NW}-\eqref{C}-\eqref{K} belong to 
$$H^3(0,T; L^{\infty}(\Lambda)) \cap H^1(0,T; H^2(\Lambda)).$$
Then there exists $T_0>0$ such that for all $T\geq T_0$, if 
$$
\ppp_x u[p] (Q,t) = \ppp_x u[q] (Q,t)  \quad \mbox{for each 
$t\in(0,T)$ and $Q\in\Pi_2\setminus\{Q_{N_2}\}$},
$$
then we have $p=q$ in $\Lambda$.
\end{theorem}

The proof of this result in Section 2 relies on a 1-d result 
of uniqueness for the determination of potential in the wave equation 
and an ``undressing'' argument.\\

It is worth mentioning that our argument gives the uniqueness for the inverse 
problems of determination of potentials on tree-shaped networks also for the 
heat and the Schr\"odinger equations using only measurements at 
$N_2 - 1$ exterior end-points. 
In fact, our arguments in proving the uniqueness for the wave and 
the Schr\"odinger equations are essentially the same and are based on local 
Carleman estimates, while the uniqueness for the inverse heat problem is
reduced to the uniqueness for the corresponding inverse wave problem (in a sense to be detailed later).

\begin{theorem}[\bf Uniqueness for the heat inverse problem]\label{Thm1-1}
Assume that $p, q \in L^{\infty}(\Lambda)$, the initial value $u^0$ satisfies 
$$ 
|u^0(x)|\geq r> 0, \quad \mbox{a.e. in $\Lambda$}
$$
for some constant $r$, and the solutions $u[p]$ and $u[q]$ to 
\eqref{NH}-\eqref{C}-\eqref{K}, belong to 
$$
H^2(0,T; L^{\infty}(\Lambda))
\cap H^1(0,T;H^2(\Lambda)).$$
Then there exists $T>0$ such that if
$$
u[p] (Q,t) = u[q] (Q,t)  \quad \mbox{for each $t\in(0,T)$ and  
$Q\in\Pi_2\setminus\{Q_{N_2}\}$},
$$
then we have $p=q$ in $\Lambda$.
\end{theorem}

\begin{theorem}[\bf Uniqueness for the Schr\"odinger inverse problem]\label{Thm1-2}
Assume that $p, q \in L^{\infty}(\Lambda)$, the initial value $u^0$ satisfies 
$$ 
|u^0(x)|\geq r> 0, \quad \mbox{a.e. in $\Lambda$}
$$
for some constant $r$, and the solutions $u[p]$ and $u[q]$ to 
\eqref{NS}-\eqref{C}-\eqref{K}, belong to 
$$H^2(0,T; L^{\infty}(\Lambda))
\cap H^1(0,T;H^2(\Lambda)).$$
Then there exists $T>0$ such that 
$$
\ppp_xu[p] (Q,t) = \ppp_xu[q] (Q,t)  \quad \mbox{for each $t\in(0,T)$ and  
$Q\in\Pi_2\setminus\{Q_{N_2}\}$},
$$
then we have $p=q$ in $\Lambda$.
\end{theorem}

One can refer to \cite{BaudouinCrepeauValein11} for the same inverse problem in the wave equation on a network where the proof is detailed in a star-shaped network but is actually generalizable to tree-shaped networks. 
Reference \cite{IgnatPazRosier12} discusses the inverse heat problem on tree-shaped
network.  Moreover the paper \cite{IgnatPazRosier12} 
treats the Schr\"odinger case in a star-shaped network and 
needs measurements at all external nodes.
We do not know any uniqueness result for non-tree graphs, which are graphs containing a closed cycle.  For observability inequality on general graph, see e.g., \cite{DagZua06}.\\

For the inverse problem in the wave equation case, we state

\begin{theorem}[\bf Stability]\label{Thm2}
Let $M>0$ and $r>0$.  Assume that $p\in L^\infty_M(\Lambda)$ and 
the solutions $u[p]$ and $u[q]$ to \eqref{NW}-\eqref{C}-\eqref{K} satisfy
$$
u[p], u[q] \in H^3(0,T;L^{\infty}(\Lambda)) \cap H^1(0,T;H^2(\Lambda)).
$$
Assume also that the initial value $u^0$ satisfies 
$$ 
|u^0(x)|\geq r> 0, \quad \mbox{a.e. in $\Lambda$}.
$$
Then there exists $T_0>0$ such that for all $T\geq T_0$, there exists 
$C=C(T,r,M, \ell_0,..., \ell_N)>0$ 
such that
\begin{equation}\label{stabpi}
|| q-p||_{L^2(\Lambda)}\leq C \sum_{j=1}^{N_2-1}\left\| \ppp_xu_{j}[p](Q_j)
- \ppp_xu_{j}[q](Q_j)\right\|_{H^1(0,T)}.
\end{equation}
\end{theorem}

This paper is composed of five sections.  The proof of uniqueness in the inverse problem in the wave equation case (Theorem~\ref{Thm1}) is presented in Section 2. The cases of Schr\"odinger and heat equations are studied in
Section 3, devoted to the proofs of Theorems~\ref{Thm1-1} and \ref{Thm1-2}.
Theorem~\ref{Thm2} is finally proven in Section 5 by a 
compactness-uniqueness argument and an observability estimate on the whole 
network.\\

We conclude this section with a classical result on the existence and regularity of solutions of the wave system and provide the corresponding energy estimates 
for the solution which we will need later.

\begin{lemma}\label{Energy}
Let $\Lambda$ be a tree-shaped network and assume that $p\in 
L^{\infty}_M(\Lambda)$, $g 
\in L^1(0,T;L^2(\Lambda))$, $u^0\in H_0^1(\Lambda)$ and $u^1\in L^2(\Lambda)$.
We consider the 1-d wave equation on the network with the conditions 
\eqref{C} and \eqref{K}:   
\begin{equation}\label{e}
\left\{\begin{array}{lll}
\ppp_t^2u-\ppp_x^2 u+p(x)u=g(x,t),&\quad\mbox{in $\Lambda\times(0,T)$},\\
		u(Q,t)=0,&\quad \mbox{in $(0,T),\, Q\in\Pi_2$},\\
u_j(P,t)=u_k(P,t),&\quad \mbox{in $(0,T),\, P\in\Pi_1,\, j,k\in
\{ 1, ..., n_P \}$}, \\
\left[ \ppp_x u (t)\right]_P=0,&\quad  \mbox{in $(0,T)$, $P\in\Pi_1$},\\
u(0) = u^0, \quad \partial_t u(0)=u^1,&\quad \mbox{in $\Lambda$}.
\end{array}\right.
\end{equation}
The Cauchy problem is well-posed and equation \eqref{e} admits a unique weak 
solution 
$$
u\in C([0,T],H_0^1(\Lambda)) \cap C^1([0,T],L^2(\Lambda)). 
$$
Moreover there exists a constant 
$C=C(\Lambda,T,M)>0$ such that for all $t\in(0,T)$, the energy 
$$
E(t) = ||\ppp_t u(t)||^2_{L^2(\Lambda)}+||\ppp_x u(t)||^2_{L^2(\Lambda)}
$$ 
of the system \eqref{e} satisfies
\begin{equation}\label{estimeenergy}
E(t) \leq C\left(||u^0||^2_{H_0^1(\Lambda)} + ||u^1||^2_{L^2(\Lambda)} +
	||g||^2_{L^1(0,T,L^2(\Lambda))} \right)
\end{equation}
and we also have the following trace estimate
\begin{equation}\label{hiddenregularity}
\sum_{j=1}^{N_2}\left\|\partial_x u_j(Q_j)\right\|_{L^2(0,T)}^2
\leq C \left(||u^0||^2_{H_0^1(\Lambda)} 
	+ ||u^1||^2_{L^2(\Lambda)} + ||g||^2_{L^1(0,T,L^2(\Lambda))} \right).
\end{equation}
\end{lemma}

The proof of the unique existence of solution to equation~\eqref{e} 
can be read in \cite[Chap. 3]{LionsMagenesBook}.  
Estimate \eqref{estimeenergy} is a classical result which can be formally 
obtained by multiplying the main equation in \eqref{e} by $\ppp_t u_{j}$, 
summing up for $j\in\left\{0,...,N\right\}$ the integral of this equality on 
$(0,T)\times e_j$  and using some integrations by parts. 
Estimate~\eqref{hiddenregularity} is a hidden regularity result which can be 
obtained by multipliers technique (we refer to \cite[Chapter 1]{Lions}). 
Formally, for the particular case of a star-shaped network of vertex $P=0$
for example, it comes from the multiplication of \eqref{e} 
by $m(x)\ppp_x u_{j}$, where $m \in C^1(\bar{\Lambda})$ with $m(0)=0$ and 
$m_j(l_j) = 1$, summing up the integrals of this equality on $(0,T)\times 
(0,l_j)$ over $j\in \{0,...,N\}$ and using integrations by parts.

\section{Uniqueness of the inverse problem - wave network case}

As already evoked in the introduction, the proof of Theorem~\ref{Thm1} will 
use a well-known 1-d result of uniqueness for the inverse problem.  We recall it in the following 
lemma.  

\begin{lemma}\label{1d}
Let $r>0$, $p\in  L^\infty(0,\ell)$ and $T > 2\ell$. Consider the 1-d wave equation in $[0,\ell]$ with homogeneous Dirichlet boundary data as follows:
\begin{equation}\label{1D}
\left\{
\begin{array}{lll}
\ppp_t^2 y - \ppp_x^2 y + p(x)y = f(x)R(x,t), \qquad& (x,t) \in (0,\ell) \times (0,T),\\
y(\ell,t)=0,&   t\in(0,T),\\
y(x,0) = 0, \ppp_{t}y(x,0)=0 ,&  x\in(0,\ell),
\end{array}
\right.
\end{equation}
where $f\in L^2(0,\ell)$ and 
$R\in H^1(0,T;L^\infty(0,\ell))$ satisfies 
$ |R(x,0)|\geq r> 0$ a.e. in $(0,\ell)$.\\
If $\, \ppp_x y (\ell,t) = 0 \,$  for all $t\in(0,T)$,
then we have $f\equiv 0$ in $(0,\ell)$ and $y\equiv 0$ in 
$(0,\ell) \times (0,T)$.
\end{lemma}

This lemma is a classical uniqueness result for the inverse source problem in a wave equation and the proof can be 
done by the method in \cite{BuKli81} on the basis of a 1-d Carleman estimate and the even extension of $y$ to negative times $t$.  We further  refer to 
Imanuvilov and Yamamoto  \cite{ImYamIP01}, \cite{ImYamCom01}, 
Klibanov \cite{Klibanov92},
Klibanov and Timonov \cite{KliTiBook} for example, and we omit details of the 
proof. \\

\noindent {\bf Proof of Theorem~\ref{Thm1}.}
We define the following operation of ``removing'' segments from the 
tree-shaped network $\Lambda$, starting from all the external nodes where we make measurements, except one. 
We divide the proof into several steps.
\\

{\bf Step 1.}  From Lemma~\ref{1d}, we can easily prove that if $e_j$ is a 
segment of $\Lambda$ which ends at an external node $Q_j\in \Pi_2$, and if 
the solutions $u[p]$ and $u[q]$ to \eqref{NW} satisfy $\ppp_x u[p](Q_j,t) 
= \ppp_x u[q](Q_j,t)$ for all  $t\in(0,T)$, then $p= q$ on the segment~$e_j$ and $u[p](x,t) 
= u[q](x,t)$ for all $x\in e_j$ and for all  $ t\in(0,T)$.
Indeed, if we set $y = u_j[p_j] - u_j[q_j]$, then 
\begin{equation}\label{1Dj}
\left\{
\begin{array}{lll}
\ppp_t^2 y - \ppp_x^2 y + p_j(x)y = (q_j-p_j)(x)u_j[q_j](x,t)\qquad &\quad (x,t)\in (0,\ell) \times (0,T),\\
y(Q_j,t)=0,&\quad   t\in(0,T),\\
y(x,0) = 0, \ppp_{t}y(x,0)=0 ,&\quad x\in(0,\ell),
\end{array}
\right.
\end{equation}
and noting that $T>0$ is sufficiently large, we can apply Lemma~\ref{1d} 
since $\ppp_x y (Q_j,t) = 0$ for 
all $t\in(0,T)$, $u_j[q_j] \in H^1(0,T;L^{\infty}(\Lambda))$ and $|u_j^0(x)|
\geq r> 0$ on $e_j$.
We obtain that $p_j \equiv q_j$ on $e_j$ and consequently $u_j[p_j](x,t) 
= u_j[p_j](x,t)$ in $e_j\times(0,T_1)$, where $T_1 \in (0,T)$ is some
constant.
\\
Therefore, for any segment $e$ with the end-points $P$ and $Q$ such that 
$Q \in \Pi_2 \setminus \{ Q_{N_2}\}$, we see that 
$p=q$ on $e$ and
$(u[p]\vert_e)(P,t) = (u[q]\vert_e)(P,t)$, $(\partial_xu[p]\vert_e)(P,t)
= (\partial_xu[q]\vert_e)(P,t)$ for $0 < t < T_1$.
Let $\Pi^2_1$ be all the interior node points $P$ of segments of
$\Lambda$ having their other end-point in $\Pi_2 \setminus \{Q_{N_2}\}$.   We note that 
$\Pi_1^2 \subset \Pi_1$.
Applying the above argument to all the exterior end-points except for
$Q_{N_2}$, we have
$$
u[p]_j(P,t) = u[q]_j(P,t), \quad 
\partial_xu[p]_j(P,t) = \partial_xu[q]_j(P,t)
$$
for each $P \in \Pi_1^2$, $0 < t < T_1$ and $j \in \{1, ..., N_3\}$.
Here by $e_1, ..., e_{N_3}$, we enumerate the segments connecting a point in $\Pi_1^2$ and a point in 
$\Pi_2 \setminus \{Q_{N_2}\}$. 
\\
{\bf Step 2.}  Let $P \in \Pi_1$ be a given node such that $n_P$ segments, say,
$e_1, ..., e_{n_P}$ meet at $P$ and $e_1, ..., e_{n_P-1}$ 
connect $P$ with exterior end-points, say, $Q_1, ..., Q_{n_P-1} \in \Pi_2$ and 
\begin{equation}\label{star}
\begin{array}{c}
u[p]_j(P,t) = u[q]_j(P,t),  \\
\ppp_x u[p]_j(P,t) = \ppp_x u[q]_j(P,t), \quad 
j \in \{1, ..., n_P-1\}, \thinspace 0<t<T.
\end{array}
\end{equation}
Using the continuity \eqref{C} and the Kirchhoff law \eqref{K} at node $P$, 
we can deduce that 
$$
\begin{array}{c}
u[p]_{n_P}(P,t) = u[q]_{n_P}(P,t),  \\
\ppp_x u[p]_{n_P}(P,t) = \ppp_x u[q]_{n_P}(P,t), \quad 0 < t < T.
\end{array}
$$

{\bf Step 3.} 
Let $\Lambda^2$ be the graph generated from $\Lambda$ by removing 
$e_1, ..., e_{N_3}$.
Therefore, since $T_1>0$ is still sufficiently large, we can 
apply the same argument as in Step 1 to the graph $\Lambda^2$.

We repeat this operation to obtain the sets $\Lambda^3$, then 
$\Lambda^4$,..., $\Lambda^n$. Hence, 
let  $L^k$ be the set of all the open segments of $\Lambda_k$,
$\Pi^k_1$ the set of the interior node points of $\Lambda_k$, $\Pi^k_2$ 
the set of external endpoints of $\Lambda_k$.
Setting $\Lambda^1 = \Lambda$, we note that $L^1 = \{ e_0, ..., e_{N}
\}$, $\Pi^1_1 = \{ P_1, ..., P_{N_1}\}$, $\Pi^1_2 = \{ Q_1, ..., Q_{N_2}\}$.

By \eqref{C} and \eqref{K}, we see that
$$
\Pi^{k-1}_1 \supset \Pi^k_1, \qquad \forall k \in \mathbb N  
$$
and 
$$
\Lambda_k = L^k \cup \Pi^k_1 \cup \Pi^k_2, \quad
L^k \cap \Pi^k_1 = L^k \cap \Pi^k_2 = \Pi^k_1 \cap \Pi^k_2 = 
\emptyset, \quad \forall k\in  \mathbb N.                           
$$
In order ro complete the proof, it is sufficient to prove there exists $n \in 
\mathbb N$ such that
\begin{equation}\label{emptylambda}
\Lambda_n = \emptyset.                
\end{equation}
Assume contrarily that $\Lambda_n \ne \emptyset$ for all $n \in \mathbb N$.  
Since every segment with exterior end-point in $\Pi_2 \setminus 
\{ Q_{N_2}\}$, can be removed (meaning that 
$u[p] = u[q]$ on the segment) by the 
above operation, we obtain that there exists $n_0 \in N$ such that 
$\Lambda_{n_0} = L^{n_0} \cup \Pi_1^{n_0}$, \textit{i.e.,} 
$\Pi^{n_0}_2 = \emptyset$.
Then $\Lambda_{n_0}$ must be a closed cycle since it possesses no external 
endpoint.
By assumption, there exist no closed cycles in a tree-shape network. 
This is a contradiction and thus the proof of \eqref{emptylambda}, and 
therefore, the one of Theorem~\ref{Thm1} is completed.
\endproof

\section{Uniqueness for the inverse problem - Schr\"odinger and heat network cases}

\subsection{Proof of Theorem \ref{Thm1-1} - Heat case.}

We apply an argument similar to the proof of Theorem 4.7 in \cite{Klibanov92} which is based on the reduction of the inverse 
heat problem to an inverse wave problem by a kind of Laplace transform
called the Reznitzkaya transform (e.g., \cite{Isakov}, 
\cite{LRS}, \cite{RomBook}).

First we define an operator $\Delta_{\Lambda}$ in $L^2(\Lambda)$ by $\Delta_{\Lambda}u = \ppp_x^2 u_j$ in $e_j$, for all $j\in \{0, 1, ...., N\}$ with
$$ 
\mathcal{D}(\Delta_{\Lambda}) = \big\{ u = (u_0, ..., u_N); \,
\mbox{$u_j \in H^2(e_j)$, $\ppp_xu(Q) = 0$ for $Q \in \Pi_2$, $u_j$ satisfying \eqref{C} and \eqref{K}}\big\}.
$$
Here, $e_j$ is oriented from $I(e_j)$ to $T(e_j)$ when defining 
$\ppp_x^2$.  Then, similarly to \cite{IgnatPazRosier12}, we can prove that $\Delta_{\Lambda}$ is self-adjoint
and $(\Delta_{\Lambda}u, u)_{L^2(\Lambda)} :=
\sum_{j=0}^N (\ppp_x^2u_j,u_j)_{L^2(e_j)} \ge 0$.  Therefore $\Delta_{\Lambda}$
generates an analytic semigroup $e^{t\Delta_{\Lambda}}$, $t>0$ (e.g.,
Pazy \cite{Pazy}, Tanabe \cite{TanabeBook}).  
Since $p \in L^{\infty}(\Lambda)$, the
perturbed operator $\Delta_{\Lambda} + p$ generates an analytic semigroup
(e.g., Theorem 2.1 in \cite{Pazy}, p.80).  Therefore by the semigroup theory 
(e.g. \cite{Pazy}, \cite{TanabeBook}), we know that the solutions $u[p](x,t)$ 
and $u[q](x,t)$ of equation  (4) are analytic in $t$ for any fixed 
$x\in \Lambda$.   More precisely, $u[p], u[q]: (0,\infty) \longrightarrow H^2(\Lambda)$ are analytic in $t>0$.

By $u^H[p]$ we denote the solution of the heat system \eqref{NH} and 
by $u^H[q]$ the corresponding solution when the potential is $q$.
By the analyticity in $t$ and the assumption in the theorem, we have
\begin{equation}\label{assumHeat}
u^H[p](Q,t) = u^H[q](Q,t), \quad \forall Q \in\Pi_2\setminus \{Q_{N_2}\}, 
\thinspace \forall t > 0.
\end{equation}
On the other hand, denote by $\widetilde u[p]$ the solution of the 
wave system 
\begin{equation}\label{NWH}
\left\{
\begin{array}{lll}
\ppp_t^2 u_j - \ppp_x^2 u_j + p_j(x)u_j = 0,\quad &\quad\forall j\in\{0,1,..., N\}, \forall (x,t)\in e_j \times (0,\infty),\\
\ppp_x u[p](Q,t)=0,&\quad \forall Q\in \Pi_2, \forall t\in(0,\infty),\\
u[p](x,0) = 0, \thinspace \ppp_t u(x,0)=u^0(x),&\quad \forall x\in\Lambda
\end{array}
\right.
\end{equation}
and by $\widetilde u[q]$ the corresponding solution when the potential is $q$.
Then we obtain (e.g., \cite[pp.251-252]{LRS}) that 
$$
\frac{1}{2\sqrt{\pi t^3}}\int^{\infty}_0 \tau e^{-\frac{\tau^2}{4t}}
\widetilde u[p](x,\tau) d\tau
$$
satisfies (4). The uniqueness of solution to equation (4) implies
$$
u^H[p](x,t) = \frac{1}{2\sqrt{\pi t^3}}\int^{\infty}_0 \tau 
e^{-\frac{\tau^2}{4t}}\widetilde u[p](x,\tau) d\tau,
\quad  \forall x \in \Lambda, \forall t > 0
$$
and the same equality with $q$.
By assumption \eqref{assumHeat}, we obtain
$$
\frac{1}{2\sqrt{\pi t^3}}\int^{\infty}_0 \tau 
e^{-\frac{\tau^2}{4t}}(\widetilde u[p] - \widetilde u[q])(Q,\tau) d\tau = 0, 
\quad
\forall Q\in \Pi_2 \setminus \{Q_{N_2}\}, \forall t>0.
$$
By the change of variables $s = \frac{1}{4t}$ and $\tau^2 = \eta$, we obtain
$$
\int^{\infty}_0 e^{-s\eta}(\www{u}[p] - \www{u}[q])(Q,\sqrt{\eta}) d\eta = 0, 
\quad \forall Q\in \Pi_2 \setminus \{Q_{N_2}\}, \forall s>0
$$
and the injectivity of the Laplace transform yields
\begin{equation}
(\www{u}[p] - \www{u}[q])(Q,\sqrt\eta) = 0, \quad \forall Q\in \Pi_2 \setminus 
\{Q_{N_2}\}, \forall \eta>0.
\end{equation}
Applying the same argument as in Section 2 for the wave system, we prove $p=q$ 
in $\Lambda$.  Thus the proof of Theorem~\ref{Thm1-1} is completed.
%
\subsection{Proof of Theorem~\ref{Thm1-2} - Schr\"odinger case.}

It is sufficient to prove the following lemma.
\\
\begin{lemma}\label{schrodi}
Let $r>0$ and $p \in L^{\infty}(0,\ell)$, $f\in L^2(0,\ell)$ be 
real-valued, and $T>0$ be arbitrarily fixed.  
We consider a 1-d Schr\"odinger equation:
$$
\left\{
\begin{array}{lll}
i\ppp_t y - \ppp_x^2 y + p(x)y = f(x)R(x,t), \qquad 
&\forall (x,t) \in (0,\ell) \times (0,T), \\
y(\ell,t) = 0, &\forall t \in (0, T),\\
y(x,0) = 0, &\forall  x \in (0, \ell),
\end{array}
\right.
$$
where $R \in H^1(0,T; L^{\infty}(0,\ell))$ satisfies $\vert R(x,0)\vert
\ge r > 0$ a.e. in $(0,\ell)$. \\
 If $\, \ppp_xy(\ell,t) = 0\,$ for all 
$ t \in (0, T)$,
then we have $f=0$ in $(0,\ell)$ and $y=0$ in $(0,\ell) \times (0,T)$.
\end{lemma}

Using the same method as the one for the proof of Lemma~\ref{1d}, this lemma is proved by means of 
the following Carleman estimate:
\\
\begin{lemma}
For $x_0 \not\in [0,\ell]$ and $\beta > 0$ arbitrarily fixed, we set 
$$
Sv = i\ppp_tv - \ppp_x^2v, \quad 
\varphi(x,t) = e^{\gamma(\vert x-x_0\vert^2 - \beta t^2)}, \quad
(x,t) \in (0,\ell) \times (0,T).
$$ 
Then there exists a constant 
$\gamma_0 > 0$ such that for arbitrary $\gamma\ge \gamma_0$ we can 
choose $s_0 > 0$ satisfying, for a constant $C>0$, 
$$
\int^T_0\int^{\ell}_0 (s\vert \ppp_xv\vert^2 + s^3\vert v\vert^2)
e^{2s\va} dxdt 
\le C\int^T_0\int^{\ell}_0 \vert Sv\vert^2 e^{2s\va} dxdt
$$
for all $s > s_0$ and all $v \in L^2(0,T;H^2_0(0,\ell)) \cap
H^1_0(0,T;L^2(0,\ell))$.
\end{lemma}

This is a Carleman estimate with regular weight function 
$\gamma(\vert x-x_0\vert^2 - \beta t^2)$ and for the proof, 
we refer to e.g. \cite[Lemma 2.1]{YuYam-AA07}
(see also \cite{YuYam-CAM10}).
Concerning a Carleman estimate for Schr\"odinger equation in a bounded domain $\Omega \subset \mathbb R^n$ with singular
weight function $\varphi$, we can refer for example to \cite{BaudouinPuel02,MOR08}.

On the basis of this lemma, the proof of Lemma~\ref{schrodi} is done by 
a usual method by Bukhgeim and Klibanov \cite{BuKli81} by 
using the extension of $y$ to $-T<t<0$ by $y(\cdot,t) 
= \overline{y(\cdot,-t)}$ and a cut-off argument.
We omit the details of the proof.

\section{Observability in the wave network}

The proof of the stability result will rely strongly on the classical 
result of 
observability that we are now presenting and proving.
One should specifically mention 
the survey \cite{ZuaSurveyNetworks} and the books \cite{DagZua06},
\cite{LaLeSch94}, where the question of observability in networks of strings (or wave equations) is widely explored in different cases.

We concentrate here on the case where the observation available comes from all but one external nodes, in a setting with a system of wave equations with potential. Since most of the literature on string networks focus only on the wave equation without lower order terms (see \cite{LaLeSch94} or \cite{DagZua06} for instance), we detail here how to obtain the observability result for the wave equation with potential.
In some other cases, we can prove the observability inequality directly
by a global Carleman estimate (e.g. \cite{BaudouinCrepeauValein11}). 

\begin{theorem}[\bf Observability inequality]\label{Thm3}
On the tree-shaped network  $\Lambda$, assuming $p\in L^\infty(\Lambda)$, let us consider the system of 1-d wave 
equations under the continuity and Kirchhoff law's assumptions
\eqref{C} and \eqref{K}:
\begin{equation}\label{eqobs}
	\left\{\begin{array}{lll}
\ppp_t^2u-\ppp_x^2 u+p(x)u=0,&\quad\mbox{in $\Lambda\times(0,T)$},\\
		u(Q,t)=0,&\quad \mbox{in $(0,T),\, \forall Q\in\Pi_2$},\\
u_j(P,t)=u_k(P,t),&\quad \mbox{in $(0,T),\,\forall P\in\Pi_1,\, \forall 
j,k\in \{1, ..., n_P\}$}, \\
\left[ \ppp_x u (t)\right]_P=0,&\quad  \mbox{in $(0,T), \forall P\in\Pi_1$},\\
u(x,0) = 0, \quad \partial_t u(x,0)=a(x),&\quad \mbox{in $\Lambda$},
\end{array}\right.
\end{equation}
Then there exists a minimal time $T_0$ such that for all $T > T_0$,
the observability 
estimate
\begin{equation}\label{obs}
\int_\Lambda | a(x)|^2 dx
\le C\sum_{j=1}^{N_2-1} \int^T_0 | \ppp_xu_j(Q_j,t)|^2 dt
\end{equation}
holds for a solution $u$ of \eqref{eqobs}.
\end{theorem}

\noindent {\bf Proof of Theorem~\ref{Thm3}.}
Let $v$ be the solution of the system
$$
\left\{
\begin{array}{lll}
\ppp_t^2 v - \ppp_x^2 v=  -pu \quad &\quad\forall (x,t)\in \Lambda \times 
(0,T),\\
v(Q,t)=0,&\quad \forall Q\in \Pi_2, t\in(0,T),\\
v_j(x,0) = 0, \ppp_t v_j(x,0)=0,&\quad\forall j\in\{0,1,..., N\},
\quad x\in e_j,
\end{array}
\right.
$$
under conditions \eqref{C} and \eqref{K}.
Then \eqref{hiddenregularity}  in Lemma~\ref{Energy} and 
$p\in L^\infty(\Lambda)$ yields
\begin{equation}\label{eqeq}
\sum_{j=1}^{N_2} \int^T_0 | \ppp_xv_j(Q_j,t)|^2dt 
\le C \int^T_0\int_\Lambda  | pu|^2 dxdt 
\le C \int^T_0\int_\Lambda | u|^2 dxdt.
\end{equation}
Setting $w = u - v$, we still have \eqref{C} and \eqref{K} satisfied by $w$, along with the following equation
$$
\left\{
\begin{array}{lll}
\ppp_t^2 w - \ppp_x^2 w=  0 \quad &\quad\forall (x,t)\in \Lambda \times 
(0,T),\\
w(Q,t)=0,&\quad \forall Q\in \Pi_2, t\in(0,T),\\
w_j(x,0) = 0, \ppp_t w_j(x,0)=a(x),&\quad\forall x\in \Lambda.
\end{array}
\right.
$$
Therefore, using a classical observability inequality in the case where $p=0$ 
(e.g., \cite{DagZua06,LaLeSch94}), we have
$$
\int_\Lambda  | a(x)|^2  dx
\le C \sum_{j=1}^{N_2-1} \int^T_0 | \ppp_xw_j(Q_j,t)|^2 dt.
$$
Hence, by \eqref{eqeq}, we have 
\begin{eqnarray}
&& \int_\Lambda | a(x)|^2 dx
\le C\sum_{j=1}^{N_2-1} \int^T_0 | \ppp_xu_j(Q_j,t)|^2 dt
+ C\sum_{j=1}^{N_2-1} \int^T_0 | \ppp_xv_j(Q_j,t)|^2 dt\nonumber\\
&&\le C\sum_{j=1}^{N_2-1} \int^T_0 | \ppp_xu_j(Q_j,t)|^2 dt
+ C \int^T_0\int_\Lambda  | u|^2dxdt.\label{unicompact}
\end{eqnarray}
Therefore a usual compactness-uniqueness argument yields the observability 
inequality \eqref{obs}.
Indeed, if \eqref{obs} is not satisfied, then we can assume 
that there exists $a^n\in L^2(\Lambda), n\in \mathbb N$ such that
\begin{equation}\label{hypcontr}
\|a^n\|_{L^2(\Lambda)} = 1, \, \forall n\in\mathbb N \quad  \hbox{ and }  \quad
\lim_{n\to + \infty}\sum_{j=1}^{N_2-1} \int^T_0 | \ppp_xu_j^n(Q_j,t)|^2 dt = 0.
\end{equation}
Using the energy estimate \eqref{estimeenergy} of Lemma~\ref{Energy} on the 
solution $u^n$ of  system \eqref{eqobs} with initial data $a^n$, we obtain
$$
|| u^n(t)||^2_{H^1_0(\Lambda)}  = ||\ppp_x u^n(t)||^2_{L^2(\Lambda)} 
\leq C||a^n||^2_{L^2(\Lambda)} \leq C.
$$
Since the embedding $H^1_0(\Lambda) \subset L^2(\Lambda)$ is compact, we can 
extract a subsequence, denoted again by the same notation and we have 
$(u^n)_{n\in\mathbb N^*}$ convergent in $L^2(\Lambda)$. Therefore, using 
\eqref{unicompact}, we obtain 
\begin{align*}
\int_\Lambda | a^n-a^m|^2 dx
\le & ~C\sum_{j=1}^{N_2-1} \int^T_0 | \ppp_xu_j^n(Q_j,t) |^2 dt
+ C\sum_{j=1}^{N_2-1} \int^T_0 |\ppp_xu_{j}^m(Q_j,t)|^2 dt \\
&
+~ C \int^T_0\int_\Lambda  | u^n - u^m|^2dxdt
\end{align*}
so that \eqref{hypcontr} and 
$
\displaystyle\lim_{n,m\to \infty} \|u^n - u^m\|_{L^2(\Lambda)} = 0
$ imply 
$
\displaystyle\lim_{n,m\to \infty} ||a^n-a^m||^2_{L^2(\Lambda)} = 0.
$
Consequently,  there exists a limit $a_0$ such that 
$\displaystyle\lim_{n\to +\infty}a^n = a_0$ in $L^2(\Lambda)$ and from 
\eqref{hypcontr} , we have $\|a_0\|_{L^2(\Lambda)} = 1$. 
Moreover, the solution $u[a_0]$ of system  \eqref{eqobs} with initial data 
$a_0$ is such that 
$$
\ppp_xu_{j}^m[a_0](Q,t) = 0 , \quad \forall t\in(0,T), \forall Q\in\Pi_2.
$$
Hence we apply a classical unique continuation result for a wave equation to obtain that
$u[a_0]$ vanishes everywhere so that $a_0 = 0$, which contradicts 
$\|a_0\|_{L^2(\Lambda)} = 1$.  Here, the unique continuation can be proved for instance 
by a Carleman estimate (e.g. \cite{Isakov}, \cite{KliTiBook}).
This ends the proof of Theorem~\ref{Thm3}.

\section{Proof of the stability for the wave network inverse problem}

This section is devoted to the proof of Theorem~\ref{Thm2}.  The proof relies 
on a compactness-uniqueness argument and the observability estimate 
(Theorem~\ref{Thm3}) on the whole network. \\

Let us denote by $u[p]$ the solution of \eqref{NW} under the assumptions 
\eqref{C} and \eqref{K}. 
Henceforth we always assume the conditions \eqref{C} and \eqref{K}.
We consider $y = \partial_t\left(u[p] - u[q] \right)$ that satisfy
\begin{equation}\label{eqy}
\left\{
\begin{array}{lll}
\ppp_t^2 y - \ppp_x^2 y + q(x)y = (q-p) \partial_t u[p]\quad 
&\quad\forall (x,t)\in \Lambda \times (-T,T),\\
y(Q,t)=0,&\quad \forall Q\in \Pi_2, t\in(0,T),\\
y(x,0) = 0, \ppp_t y(x,0)=(q-p)u^0(x),&\quad \forall x\in\Lambda,
\end{array}
\right.
\end{equation}
We define $\psi$ and $\phi$ as the solutions of
\begin{equation}\label{eqpsi}
\left\{
\begin{array}{lll}
\ppp_t^2 \psi - \ppp_x^2 \psi + q(x)\psi = (q-p) \partial_t u[p]\quad 
&\quad\forall (x,t)\in \Lambda \times (-T,T),\\
\psi(Q,t)=0,&\quad \forall Q\in \Pi_2, t\in(0,T),\\
\psi(x,0) = 0, \ppp_t \psi(x,0)=0,&\quad \forall x\in\Lambda,
\end{array}
\right.
\end{equation}
and 
\begin{equation}\label{eqphi}
\left\{
\begin{array}{lll}
\ppp_t^2 \phi - \ppp_x^2 \phi + q(x)\phi = 0\quad &\quad
\forall (x,t)\in \Lambda \times (-T,T),\\
\phi(Q,t)=0,&\quad \forall Q\in \Pi_2, t\in(0,T),\\
\phi(x,0) = 0, \ppp_t \phi(x,0)=(q-p)u^0(x),&\quad \forall x\in\Lambda.
\end{array}
\right.
\end{equation}
such that $y = \psi + \phi$.  We can apply Theorem~\ref{Thm3} to equation 
\eqref{eqphi} so that 
\begin{equation}\label{obsphi}
\int_\Lambda | (q-p)u^0|^2 dx
\le C\sum_{j=1}^{N_2-1} \int^T_0 | \ppp_x\phi_j(Q_j,t)|^2 dt.
\end{equation}
On the other hand, a regularity result of Lemma~\ref{Energy} applied 
to a time derivative of equation \eqref{eqpsi} gives 
\begin{eqnarray}
\sum_{j=1}^{N_2}\left\|\partial_x \psi_j(Q_j)\right\|_{H^1(0,T)}^2
&\leq& C \left( ||(q-p) \partial_{t}^2 u[p]||^2_{L^1(0,T,L^2(\Lambda))} 
+ ||(q-p) u^1||^2_{L^2(\Lambda)}\right)\nonumber\\
&\leq& 2CK^2  ||q-p||^2_{L^2(\Lambda)} \label{hiddenregularitypsi}
\end{eqnarray}
as soon as we have $ u[p]\in H^2(0,T,L^\infty(\Lambda))$ which yields 
$\partial_tu[p]\in C([0,T];L^\infty(\Lambda))$ so that 
$u^1\in L^\infty(\Lambda)$) with $\|u[p]\|_{H^2(0,T,L^\infty(\Lambda))}
\leq K$. The compact embedding $H^1 (0,T) \subset L^2(0,T)$ allows then to 
write that the operator $\Psi : L^2(\Lambda) \to  L^2(0,T)$ defined by 
$$
\Psi(p-q)(t) = \sum_{j=1}^{N_2} \partial_x \psi_j(Q_j,t) , \qquad 0<t<T
$$
is compact.

Therefore, since we have $|u^0(x)|\geq r> 0$ almost everywhere in $\Lambda$, 
by \eqref{obsphi} and \eqref{hiddenregularitypsi}, we obtain
\begin{eqnarray}
|| q-p||_{L^2(\Lambda)}
&\leq& C \int_\Lambda | (q-p)u^0|^2 dx~  \leq ~ C\sum_{j=1}^{N_2-1} \int^T_0 | 
\ppp_x\phi_j(Q_j,t)|^2 dt   	\nonumber\\
&\leq& C\sum_{j=1}^{N_2-1} \int^T_0 | \ppp_x y_j(Q_j,t)|^2 dt 
+ C\sum_{j=1}^{N_2} \int^T_0 | \ppp_x\psi_j(Q_j,t)|^2 dt
                          	\nonumber\\
&\leq& C\sum_{j=1}^{N_2-1} \int^T_0 | \ppp_x y_j(Q_j,t)|^2 dt 
+ C  ||\Psi(q-p)||^2_{L^2(0,T)}  		\label{estimfinale}\\
&\leq& C \sum_{j=1}^{N_2-1}\left\| \ppp_xu_{j}[p](Q_j)- \ppp_xu_{j}[q](Q_j)
\right\|_{H^1(0,T)} + 	C ||\Psi(q-p)||^2_{L^2(0,T)}.
                                     	\nonumber
\end{eqnarray}
We aim at proving that we can get rid of the second term on the right-hand 
side of the last estimate in order to obtain \eqref{stabpi}. 
Again, a compactness-uniqueness argument will be the key and it relies here on 
the compactness of $\Psi$ and the uniqueness result of Theorem~\ref{Thm1}.

Indeed, we set $f= q-p$.  We assume that 
$$
|| f||_{L^2(\Lambda)}\leq C \sum_{j=1}^{N_2-1}\left\| \ppp_x y_j(Q_j)
\right\|_{L^2(0,T)},
$$
which is equivalent to \eqref{stabpi}, does not hold.  Then one can assume 
that there exists $f^n\in L^2(\Lambda), n\in\mathbb N$ such that
\begin{equation}\label{hypcontr1}
\|f^n\|_{L^2(\Lambda)} = 1, \, \forall n\in\mathbb N 
 \quad \hbox{ and }  \quad
 \lim_{n\to + \infty} \sum_{j=1}^{N_2-1}\left\| \ppp_x y_j^n(Q_j)
\right\|_{L^2(0,T)} = 0.
\end{equation}

First, since the sequence $(f^n)_{n\in\mathbb N}$ is bounded in 
$L^2(\Lambda)$, 
we can extract a subsequence denoted again by $(f^n)_{n\in\mathbb N}$ 
such that 
it converges towards some $f^0\in L^2(\Lambda)$ weakly in $L^2(\Lambda)$. 
Since $\Psi$ is a compact operator, we obtain therefore the strong 
convergence result
\begin{equation}\label{strongcv}
\lim_{n,m\to \infty} \|\Psi(f^n) - \Psi(f^m)\|_{L^2(0,T)} = 0.
\end{equation}
Then, from \eqref{estimfinale} we can write
$$
|| f^n - f^m||_{L^2(\Lambda)}
\leq C \sum_{j=1}^{N_2-1}\left\| \ppp_x y_j^n(Q_j)\right\|_{L^2(0,T)} 
+ C \sum_{j=1}^{N_2-1}\left\| \ppp_x y_j^m(Q_j)\right\|_{L^2(0,T)}
+ C ||\Psi(f^n) - \Psi(f^m)||^2_{L^2(\Lambda)}
$$
and deduce from \eqref{hypcontr1} and \eqref{strongcv} that $\displaystyle\lim_{n,m\to \infty} 
\|f^n - f^m\|_{L^2(\Lambda)} = 0$, so that $\displaystyle\lim_{n\to \infty} 
\|f^n - f^0\|_{L^2(\Lambda)} = 0$ with 
\begin{equation}\label{f0}
\|f^0\|_{L^2(\Lambda)} = 1.
\end{equation}

Moreover, using the trace estimate \eqref{hiddenregularity} of 
Lemma~\ref{Energy} for the solution $y^n$ of  system \eqref{eqy} with initial 
data $f^n u^0$ and source term $f^n\partial_tu[p]$, we obtain 
$$
\sum_{j=1}^{N_2-1}\left\|\partial_x y_j^n(Q_j)\right\|_{L^2(0,T)}^2
\leq C\left( ||f^n u^0||^2_{L^2(\Lambda)} 
+ ||f^n\partial_tu[p]||^2_{L^1(0,T,L^2(\Lambda))} \right)
\leq 2CK^2 \|f^n\|_{L^2(\Lambda)}.
$$
Thus we can write 
$$
\lim_{n\to \infty} \sum_{j=1}^{N_2-1}\left\|\partial_x y_j^n(Q_j) 
- \partial_x y_j^0(Q_j)\right\|_{L^2(0,T)}^2
\leq 2CK^2 \lim_{n\to \infty} \|f^n - f^0\|_{L^2(\Lambda)} = 0,
$$
which, combined with \eqref{hypcontr1}, gives 
$$
\partial_x y_j^0(Q,t) = 0, \qquad \forall Q\in \Pi_2\setminus {\{Q_{N_2}\}}, \forall t\in(0,T).
$$
We finally apply Theorem~\ref{Thm1} and obtain $f^0 = 0$ in $L^2(\Lambda)$, 
which contradicts \eqref{f0}. 
Thus the proof of Theorem~\ref{Thm2} is complete.

\bibliographystyle{plain}

\providecommand{\MR}[1]{}

\end{document}